\documentclass{article}
\usepackage{latexsym}
\usepackage{amssymb}
\newtheorem{theorem}{Theorem}

\usepackage{amsfonts,amsmath,amssymb}
\usepackage{pgf}
\usepackage{tikz}

\newcommand*{\xMin}{0}%
\newcommand*{\xMax}{26}%
\newcommand*{\yMin}{0}%
\newcommand*{\yMax}{6}%

\begin{document}

\title{\bf On ${k}$-Dyck paths with a negative boundary}
\author{\Large Helmut Prodinger  \\Department of Mathematical Sciences, Mathematics Division, \\ Stellenbosch University, Private Bag X1, 7602 Matieland,\\ South Africa,
{\tt hproding@sun.ac.za}}
\date{}
\maketitle
\begin{abstract}
	Paths that consist of  up-steps of one unit and down-steps of $k$ units, being bounded below by a horizontal line $-t$, behave like $t+1$ ordered tuples of $k$-Dyck paths, provided that $t\le k$. We describe the general case, allowing $t$ also to be larger. Arguments are bijective and/or analytic.
\end{abstract}

\section{Folklore results about $k$-Dyck paths}

\begin{figure}[h]\label{fff}
	\begin{center}
		\begin{tikzpicture}[scale=0.35]
		
		%\draw[step=1.cm,black] (-0.0,-0.0) grid (20,6);
		\draw[ultra thick] (0.0,0.) to (1.,1.);
		
		\draw (1,1) .. controls (2,5) .. (3,1);
		\draw (3,1) .. controls (4,7) .. (5,1);
		\node at (5.5,1) {$\cdot$};
		\node at (6,1) {$\cdot$};
		\node at (6.5,1) {$\cdot$};
		\draw (7,1) .. controls (8,3) .. (9,1);
		\draw [ultra thick](9,1) to  (10,2);
		\draw (1+9,1+1) .. controls (2+9,6) .. (3+9,1+1);
		\draw (3+9,1+1) .. controls (4+9,4) .. (5+9,1+1);
		\node at (5.5+9,2) {$\cdot$};
		\node at (6+9,2) {$\cdot$};
		\node at (6.5+9,2) {$\cdot$};
		\draw (16,2) .. controls (17,5) .. (18,2);	
		\draw [ultra thick](18,2) to  (19,0);
		
		\draw (19,0) .. controls (20,6) .. (21,0);
		\draw (21,0) .. controls (22,3) .. (23,0);
		\node at (23.5,0) {$\cdot$};
		\node at (24,0) {$\cdot$};
		\node at (24.5,0) {$\cdot$};
		\draw (25,0) .. controls (26,4) .. (27,0);

		\end{tikzpicture}
		\end	{center}
		
		\caption{The decomposition of generalized Dyck paths for $k=2$. Reading from left to right, the decomposition leads to
			$xyxyxy$.}
		
		\end{figure}
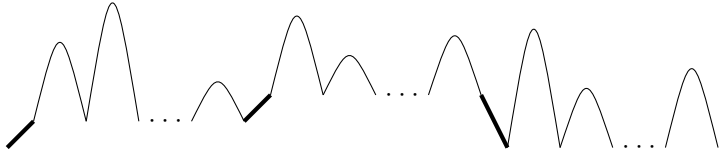

A Dyck path consists of up-steps and down-steps, one unit each, starts at the origin and returns to the 
origin after $2n$ steps, and never goes below the $x$-axis. The enumeration involves the ubiquitous Catalan numbers~\cite{Stanley-Catalan}. The family of $k$-Dyck paths is defined similarly, but the down-steps are now by $k$ units in one step. Practically every book on combinatorics has something about this; we only give two citations: \cite{BaFl, FS}. The generating function $y=y(z)=y_k(z)$ of these objects, according to length (the number of steps) can be found by a first return to the $x$-axis decomposition:
\begin{equation*}
y=1+(zy)^k\cdot z\cdot y=1+z^{k+1}y^{k+1},
\end{equation*}
where the last $z$ represents the down-step that brings the path back to the $x$-axis for the first time. The Figure \ref{fff} describes this readily; the term `1' refers to the empty path (of length 0).
(An equivalent concept is the family of $(k+1)$-ary trees; there are bijections between paths and trees, and various parameters translate accordingly.)

One can only return to the $x$-axis after a multiple of $k+1$ steps, as each down-step requires $k$ up-steps for compensation.
Consequently we may write $x=z^{k+1}$. 
Furthermore we set  $y=1+w$, making the equation  amenable to the Lagrange inversion:
\begin{equation*}
w=x(1+w)^{k+1},
\end{equation*}
and we can compute the coefficients of $w$;
\begin{equation*}
[x^n]w=\frac 1n[w^{n-1}](1+w)^{(k+1)n}=\frac 1n\binom{(k+1)n}{n-1}.
\end{equation*}
Now we compute
\begin{align*}
[x^n]y^j&=\frac1{2\pi i}\oint\frac{dx}{x^{n+1}}(1+w)^j\\
&=\frac1{2\pi i}\oint\frac{dw(1+w)^{(n+1)(k+1)}}{w^{n+1}}(1+w)^j\frac{1-kw}{(1+w)^{k+2}}\\
&=[w^n](1-kw)(1+w)^{n(k+1)-1+j}\\
&=\binom{n(k+1)-1+j}{n}-k\binom{n(k+1)-1+j}{n-1}\\
&=\frac{j}{(k+1)n+j}\binom{(k+1)n+j}{n};
\end{align*}
we will use these coefficients of powers of $y=y_k(z)$ in the next section. The contour is a small circle in the $x$-plane; the substitution from $x$ to $w$ does not change the winding number. Such computations are quite common in the context of lattice paths and/or trees.

\section{$k_t$-Dyck paths}

Selkirk~\cite{Selkirk-master} introduced an extra parameter $t$ to the family of $k$-Dyck paths. The paths might go below the $x$-axis, but never go below the horizontal line $-t$.

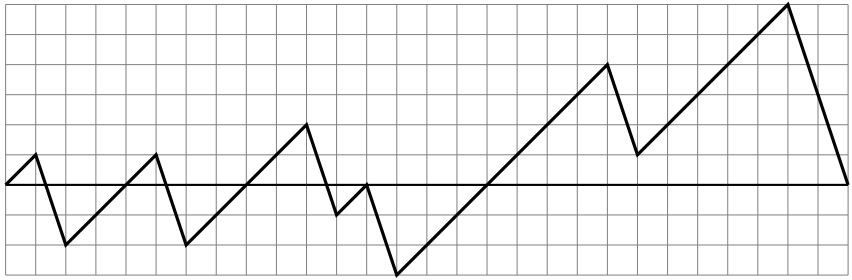
\begin{figure}[h!]
	\begin{center}
		\begin{tikzpicture}[scale=0.40]
		\foreach \i in {\xMin,...,28} {
			\draw [very thin,gray] (\i,-3) -- (\i,6)  ;
		}

		\foreach \i in {-3,...,6} {
			\draw [very thin,gray] (0,\i) -- (28,\i) ;
		}
		
		\draw [ thick] (0,0) -- (28,0)  ;
		
		\draw[very thick](0,0)--(1,1)--(2,-2)--(3,-1)--(4,0)--(5,1)--(6,-2)--(7,-1)--(8,0)--(9,1)--(10,2)--(11,-1)--(12,0)--(13,-3)--(14,-2)--(15,-1)--(16,0)--(17,1)--(18,2)--(19,3)--(20,4)--(21,1)--(22,2)--(23,3)--(24,4)--(25,5)--(26,6)--(27,3)--(28,0);		
		\end{tikzpicture}
	\end{center}
	\caption{\emph{A $3_3$-Dyck path: down-steps of 3 units and bounded below by the line $-3$.}}
	\label{path1} 
\end{figure}

The enumeration of $k_t$-Dyck paths is as follows~\cite{Selkirk-master}:

\begin{theorem}
	For $0\le t\le k$, the number of $k_t$-Dyck paths of length $(k+1)n$ is given by
	\begin{equation*}
	\frac{t+1}{(k+1)n+t+1}\binom{(k+1)n+t+1}{n}.
	\end{equation*}
	Equivalently, the generating function of $k_t$-Dyck paths by length is given by $y^{t+1}$.
\end{theorem}
That $y^{t+1}$ has indeed these coefficients was discussed in the previous section.
It also enumerates ordered $(t+1)$-tuples of $k$-Dyck paths, and the bijection in \cite{Selkirk-master} is between these two families of objects. It is to be noted that \cite{GPW} contains  somewhat equivalent statements in the language of $(k+1)$-ary trees, and instead of a boundary line, the nodes are coloured, and the colour of the root plays a role similar to the $t$ in $k_t$-Dyck paths. We report this information from \cite{Selkirk-master}.

In the present note we want to look at this again, but we also want to explain what happens if the condition
$0\le t\le k$ is no longer satisfied, i.e., if $t>k$ is allowed. The generating function is no longer $y^{t+1}$, and has to be replaced by something more complicated.

But let us start with some bijective arguments. While we consider  $3$-Dyck paths in the following figures, the illustrations are representative for other values of $k$ as well.

\begin{center}
\begin{figure}[h!]\label{figtwo}
	\begin{center}
		\begin{tikzpicture}[scale=0.365]
		\foreach \i in {\xMin,...,28} {
			\draw [very thin,gray] (\i,-3) -- (\i,6)  ;
		}
		\foreach \i in {-3,...,0} {
			\draw [very thin,green] (\i,-3) -- (\i,6)  ;
		}
		\foreach \i in {0,...,6} {
			\draw [very thin,gray] (0,\i) -- (0,\i) ;
		}
		\foreach \i in {-3,...,6} {
			\draw [very thin,green] (-3,\i) -- (0,\i) ;
		}
		\foreach \i in {-3,...,6} {
			\draw [very thin,gray] (0,\i) -- (28,\i) ;
		}
		
		\draw [ thick] (0,-3) -- (28,-3)  ;
		\draw [ thick,green] (-3,-3) -- (0,-3)  ;
		
		\draw[very thick](0,0)--(1,1)--(2,-2)--(3,-1)--(4,0)--(5,1)--(6,-2)--(7,-1)--(8,0)--(9,1)--(10,2)--(11,-1)--(12,0)--(13,-3)--(14,-2)--(15,-1)--(16,0)--(17,1)--(18,2)--(19,3)--(20,4)--(21,1)--(22,2)--(23,3)--(24,4)--(25,5)--(26,6)--(27,3)--(28,0);
		\draw[very thick,red ](-3,-3)--(0,0)		;
		\end{tikzpicture}
	\end{center}
	\caption{\emph{The   path from Figure~\ref{path1}, lifted up 3 units, with a sequence of up-steps in the beginning.}}
	\label{path2}
\end{figure}
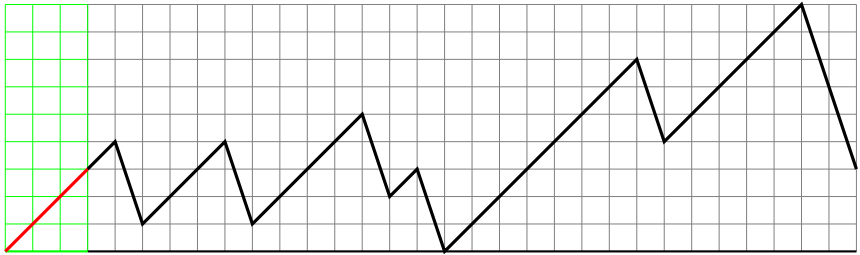
\end{center}

\begin{figure}[h!]
	\begin{center}
		\begin{tikzpicture}[scale=0.365]
		\foreach \i in {\xMin,...,28} {
			\draw [very thin,gray] (\i,-3) -- (\i,6)  ;
		}
		\foreach \i in {-3,...,0} {
			\draw [very thin,gray] (\i,-3) -- (\i,6)  ;
		}
		\foreach \i in {-3,...,6} {
			\draw [very thin,gray] (-3,\i) -- (0,\i) ;
		}
		\foreach \i in {-3,...,6} {
			\draw [very thin,gray] (0,\i) -- (28,\i) ;
		}
		
		\draw [ thick] (0,-3) -- (28,-3)  ;
		\draw [ thick] (-3,-3) -- (0,-3)  ;
		
		\draw[very thick](0,0)--(1,1)--(2,-2)--(3,-1)--(4,0)--(5,1)--(6,-2)--(7,-1)--(8,0)--(9,1)--(10,2)--(11,-1)--(12,0)--(13,-3)--(14,-2)--(15,-1)--(16,0)--(17,1)--(18,2)--(19,3)--(20,4)--(21,1)--(22,2)--(23,3)--(24,4)--(25,5)--(26,6)--(27,3)--(28,0);
		\draw[very thick](-3,-3)--(0,0)		;
		\draw[ultra thick, blue](13,-3)--(14,-2);
		\draw[ultra thick, blue](14,-2)--(15,-1);
		\draw[ultra thick, blue] (15,-1)--(16,0);
		\draw[dotted,red,very thick](14,-2)--(28,-2);
		\draw[dotted,red,very thick](14,-2)--(14,6);
		\draw[dotted,red,very thick](15,-1)--(28,-1);
		\draw[dotted,red,very thick](15,-1)--(15,6);
		\draw[dotted,red,very thick](16,0)--(28,0);
		\draw[dotted,red,very thick](16,0)--(16,6);
		\end{tikzpicture}
	\end{center}
	\caption{\emph{The path from Figure~\ref{path2}, decomposed.}}
	\label{path3}
\end{figure}
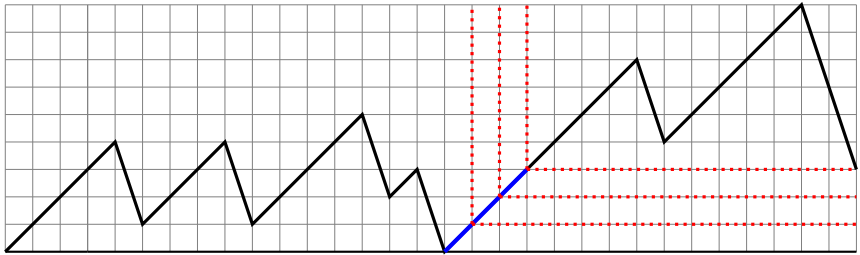

We lift up a $k_t$-Dyck paths by $t$ units and add $t$ up-steps in the beginning, as can be seen in Figure~\ref{figtwo}. The resulting path is
a $k$-Dyck path, but does not end on level 0, but rather on level $t$. It is classical (see \cite[page 321]{FS}) that these paths have generating function $z^ty^{t+1}$, thanks to a decomposition that is sketched in Figure~\ref{path3}. The first part, according to this decomposition,
ends, where the $x$-axis (=level 0) is visited for the last time. After an up-step, the second part starts and ends when the level 1 is visited for the last time, and so on. All these $t+1$ parts are $k$-Dyck paths themselves, and altogether $t$ up-steps have been identified.

Removing $t$ extra up-steps, we are at the generating function $y^{t+1}$ again, and the decomposition gives us $t+1$ (ordered) $k$-Dyck paths. In the example, we get 4 paths, see Fig.~5:
		\begin{center}
	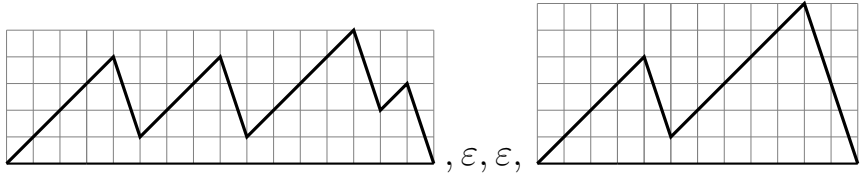
\begin{figure}[h!]	\label{path4}

			\begin{tikzpicture}[scale=0.355]
			%\hspace*{-2cm}
			\foreach \i in {\xMin,...,13} {
				\draw [very thin,gray] (\i,-3) -- (\i,2)  ;
			}
			\foreach \i in {-3,...,0} {
				\draw [very thin,gray] (\i,-3) -- (\i,2)  ;
			}
			\foreach \i in {-3,...,2} {
				\draw [very thin,gray] (-3,\i) -- (0,\i) ;
			}
			\foreach \i in {-3,...,2} {
				\draw [very thin,gray] (0,\i) -- (13,\i) ;
			}
			
			\draw [ thick] (0,-3) -- (13,-3)  ;
			\draw [ thick] (-3,-3) -- (0,-3)  ;
			
			\draw[very thick](0,0)--(1,1)--(2,-2)--(3,-1)--(4,0)--(5,1)--(6,-2)--(7,-1)--(8,0)--(9,1)--(10,2)--(11,-1)--(12,0)--(13,-3);
			\draw[very thick](-3,-3)--(0,0)		;
			\end{tikzpicture}
			\begin{tikzpicture}[scale=0.365]
			\Large$\ ,\varepsilon,\varepsilon,\ $
			\end{tikzpicture}
			\begin{tikzpicture}[scale=0.355]
			\hspace*{1.1cm}
			\begin{scope}[shift={(5,0)}]
			\foreach \i in {\xMin,...,9} {
				\draw [very thin,gray] (\i,-3) -- (\i,3)  ;
			}
			\foreach \i in {-3,...,3} {
				\draw [very thin,gray] (\i,-3) -- (\i,3)  ;
			}
			\foreach \i in {-3,...,3} {
				\draw [very thin,gray] (-3,\i) -- (0,\i) ;
			}
			\foreach \i in {-3,...,3} {
				\draw [very thin,gray] (0,\i) -- (9,\i) ;
			}
			
			\draw [ thick] (0,-3) -- (9,-3)  ;
			\draw [ thick] (-3,-3) -- (0,-3)  ;
			
			\draw[very thick](0,0)--(1,1)--(2,-2)--(3,-1)--(4,0)--(5,1)--(6,2)--(7,3)--(8,0)--(9,-3);
			\draw[very thick](-3,-3)--(0,0)		;
			\end{scope};
			\end{tikzpicture}

		\caption{Decomposed into 4 paths; the second and third paths are the empty path $\varepsilon$.}
	\end{figure}
		\end{center}	
The paper \cite{Selkirk-master} provides the same decomposition into $(t+1)$ $k$-Dyck paths. We hope that our alternative description is natural and easy to understand.

Since in this running example, two copies of the empty path appear, we provide an additional example:
\begin{figure}[h!]
	\begin{center}
		\begin{tikzpicture}[scale=0.25]
		\foreach \i in {\xMin,...,35} {
			\draw [very thin,gray] (\i,0) -- (\i,6)  ;
		}
		\foreach \i in {0,...,6} {
			\draw [very thin,gray] (0,\i) -- (35,\i)  ;
		}

		\draw [ thick] (0,0) -- (35,0)  ;

		\draw[very thick](0,0)--(1,1)--(2,2)--(3,3)--(4,0)--(5,1)--(6,2)--(7,3)--(8,4)--(9,5)--(10,2)--(11,3)--(12,4)--(13,5)--(14,2)--(15,3)--(16,4)--(17,1)--(18,2)--(19,3)--(20,4)--(21,5)--(22,2)--(23,3)--(24,4)--(25,5)--(26,6)--(27,3)--(28,4)--(29,5)--(30,6)
		--(31,3)--(32,4)--(33,5)--(34,6)--(35,3);
		
		\draw[ultra thick, blue](4,0)--(5,1);
		\draw[ultra thick, blue](17,1)--(18,2);
		\draw[ultra thick, blue] (22,2)--(23,3);
		\draw[dotted,red,very thick](5,6)--(5,1);
		\draw[dotted,red,very thick](5,1)--(35,1);
		\draw[dotted,red,very thick](18,6)--(18,2);
				\draw[dotted,red,very thick](18,2)--(35,2);
		\draw[dotted,red,very thick](23,3)--(23,6);
		\draw[dotted,red,very thick](23,3)--(35,3);
		%\draw[dotted,red,very thick](16,0)--(16,6);
				\end{tikzpicture}
				
				\begin{tikzpicture}[scale=0.25]
				\foreach \i in {\xMin,...,4} {
					\draw [very thin,gray] (\i,0) -- (\i,6)  ;
				}
				\foreach \i in {0,...,6} {
					\draw [very thin,gray] (0,\i) -- (4,\i)  ;
				}

				\draw [ thick] (0,0) -- (4,0)  ;

				\draw[very thick](0,0)--(1,1)--(2,2)--(3,3)--(4,0);
								\end{tikzpicture}\hspace*{2.7pt}
				\begin{tikzpicture}[scale=0.25]
				\foreach \i in {\xMin,...,12} {
					\draw [very thin,gray] (\i,0) -- (\i,6)  ;
				}
				\foreach \i in {0,...,6} {
					\draw [very thin,gray] (0,\i) -- (12,\i)  ;
				}

				\draw [ thick] (0,0) -- (12,0)  ;

				\draw[very thick](5-5,1-1)--(6-5,2-1)--(7-5,3-1)--(8-5,4-1)--(9-5,5-1)--(10-5,2-1)--(11-5,3-1)--(12-5,4-1)--(13-5,5-1)--(14-5,2-1)--(15-5,3-1)--(16-5,4-1)--(17-5,1-1);

				\end{tikzpicture}\hspace*{2.7pt}
				\begin{tikzpicture}[scale=0.25]
				\foreach \i in {\xMin,...,4} {
					\draw [very thin,gray] (\i,0) -- (\i,6)  ;
				}
				\foreach \i in {0,...,6} {
					\draw [very thin,gray] (0,\i) -- (4,\i)  ;
				}

				\draw [ thick] (0,0) -- (4,0)  ;

				\draw[very thick](0,0)--(1,1)--(2,2)--(3,3)--(4,0);
				\end{tikzpicture}\hspace*{2.7pt} 
				\begin{tikzpicture}[scale=0.25]
				\foreach \i in {\xMin,...,12} {
					\draw [very thin,gray] (\i,0) -- (\i,6)  ;
				}
				\foreach \i in {0,...,6} {
					\draw [very thin,gray] (0,\i) -- (12,\i)  ;
				}

				\draw [ thick] (0,0) -- (12,0)  ;

				\draw[very thick](0,0)--(1,1)--(2,2)--(3,3)--(4,0);
				\draw[very thick](0+4,0)--(1+4,1)--(2+4,2)--(3+4,3)--(4+4,0);
				\draw[very thick](0+4+4,0)--(1+4+4,1)--(2+4+4,2)--(3+4+4,3)--(4+4+4,0);
				\end{tikzpicture}
	\end{center}
	\caption{\emph{A path for $k=3$ and its decomposition.}}
	\label{path3a}
\end{figure}
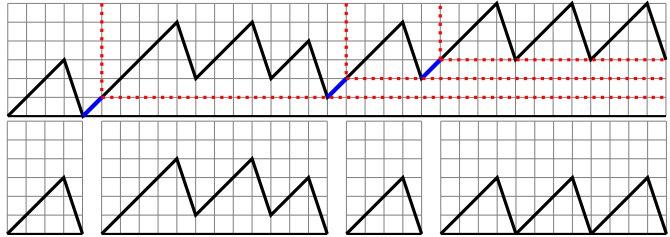

\vbox{
Observe that the operation ``\emph{shifting up the path by $t$ units}'' makes the \emph{old} origin the \emph{first point where the level $t$ is reached}. In the beginning, there are the $t$ extra up steps; the \emph{new} origin is at the left end of these added steps. If $t\le k$ this is indeed the only option to reach this level for the first time without going below the $x$-axis by using a down-step.

If $t$ is arbitrarily large, this is no longer true: We can ``live'' in the strip with boundaries 0 and $t-1$, end at the highest level, and make one up-step to reach the $k_t$-path. Fig.~\ref{path5} is a   drawing explaining this.}

\begin{center}
\renewcommand*{\yMax}{7}%
\begin{figure}[h!]
	
		\begin{tikzpicture}[scale=0.40]
		\foreach \i in {\xMin,...,\xMax} {
			\draw [very thin,gray] (\i,\yMin) -- (\i,\yMax)  ;
		}

		\foreach \i in {\yMin,...,\yMax} {
			\draw [very thin,gray] (\xMin,\i) -- (\xMax,\i) ;
		}
		
		\draw [ thick] (0,4) -- (\xMax,4)  ;
		
		\path	 node(g2) at (0,4)  [left]{$t$};
		\path	 node(g3) at (0,3)  [left]{$t-1$};
		\path	 node(g4) at (8,0)  [below]{$j$};
		\path	 node(g5) at (26,0)  [below]{$n$};
		\path	 node(g6) at (0,0)  [below]{$0$};
		
		\draw[thick](0,0)--(2,2)--(3,1) --(5,3)--(6,2)--(7,3);
		\draw[ultra thick,red](7,3)--(8,4);
		\draw (8,4) .. controls (12,8) .. (16,2) .. controls (18,0) and (20,8)   .. (26,4);
		
		\draw[<-,thick] (0,-1)-- (4,-1);	\draw[->,thick]  (4,-1)--(8,-1);
		\path	 node(g5) at (4,-1)  [below]{$F$};
		
		\draw[<-,thick] (8,-1)-- (17,-1);	\draw[->,thick]  (17,-1)--(26,-1);
		\path	 node(g6) at (17,-1)  [below]{$G$};
		\end{tikzpicture}
\label{FandG}
	\caption{The path has a first part $F$ and a second part $G$. The length of $F$ will later be called $\mathcal{J}$ and analyzed for $k=1$.}
	\label{path5}
\end{figure}
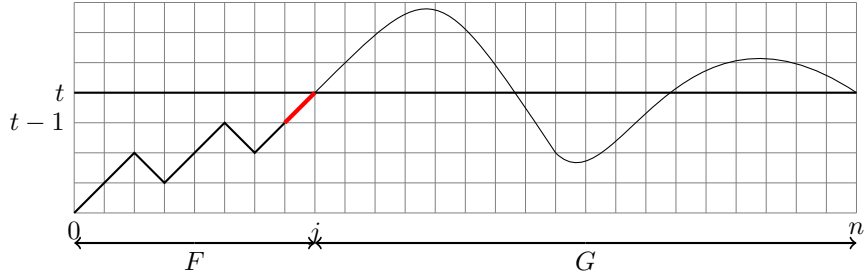
	\end{center}
	
	The decomposition of $k$-Dyck paths ending on level $t$ into $F$ and $G$ (as shown in Figure~\ref{path5}) is canonical: $G$ starts when the level $t$ has been reached for the first time, and what comes before is called $F$; this $F$ ends with an up-step, and the part before ``lives'' in the strip $0..t-1$.
	
So it is evident that for general $t$ the first part $F$ has more freedom. In the next section we will use generating functions to understand the roles of $F$ and $G$ better. In particular, the length of the $F$ may vary now.

\section{Generating functions}

At the moment, we leave the bijective arguments and go back to the original question about $k$-Dyck paths with the negative boundary $-t$. In order to obtain the relevant generating functions, we also introduce (temporarily) an upper boundary at $h$.
This means that we consider paths, living in the strip $-t..h$. 
This has the advantage that the generating functions that will appear are rational. We consider generating functions $\varphi_i(z)$, where $i$ marks the level of the endpoint;
\begin{multline*}
\varphi_i(z)=\sum_{n\ge0}z^n[\text{number of $k$-Dyck paths of length $n$,}\\[-0.5cm]
\text{ bounded by $-t$ and $h$, ending on level $i$}].
\end{multline*}

The recursion $\varphi_i=z\varphi_{i-1}+z\varphi_{i+k}+[\![i=0]\!]$, provided all indices are within the interval $-t..h$, is easy to understand; if one ends on level $i$, one must have been at level $i-1$ or at level $i+k$ before the last step. This is best written as a linear system
\begin{equation*}
\left[\begin{matrix}
1 & 0&0&\dots& 0&-z&\dots \\
-z&1 & 0&0&\dots& 0&-z&\dots \\
0&-z&1 & 0&0&\dots& 0&-z&\dots \\
0&0&-z&1 & 0&0&\dots& 0&-z&\dots \\
&&&&&\vdots\\&&&&&\vdots\\
0&0&0&0&0&0&\dots& 0&-z&1
\end{matrix}\right]
\left[\begin{matrix}
\varphi_{-t}\\
\varphi_{-t+1}\\
\varphi_{-t+2}\\
\vdots\\
\varphi_{h-2}\\
\varphi_{h-1}\\
\varphi_{h}\\
\end{matrix}\right]
=
\left[\begin{matrix}
0\\
0\\
0\\
1\\
0\\
0\\
0\\
\end{matrix}\right]
\end{equation*}

The entry 1 on the righthand side corresponds to the function $\varphi_0$. Let $D_m$ be the determinant of the matrix with $m$ rows (and columns). We find by  expanding along the first column, say, the recursion
\begin{equation*}
D_m=D_{m-1}-z^{k+1}D_{m-k-1},\quad D_0=D_1=\dots=D_{k}=1.
\end{equation*} 
The solution is
\begin{equation*}
D_m=\sum_{0\le \ell \le \frac mk}\binom{m-k\ell}{\ell}z^{(k+1)\ell}(-1)^\ell.
\end{equation*}
This can be proved  by induction on $m$ or by other methods (possibly by the use of Lambert's and Lagrange's trinomial equation; see \cite{GKP}, but this has not been investigated). For $k=1$, these polynomials are sometimes called Fibonacci polynomials and appear e.~g. in \cite{deBrKnRi72}; this highly cited paper deals with the height of planar trees, and this is equivalent to the height of Dyck paths. For $k=2$, they appear in \cite{Naiomi-paper} when computing generating functions related to $2$-Dyck paths with a boundary.
 A general reference about this method is \cite{Prodinger-handbook}.

We want to link the polynomials $D_m$ to the generating function $y=y_k$: The characteristic equation for the sequence
$D_m$ is $\lambda^{k+1}=\lambda^k-z^{k+1}$. The equation satisfied by $y$ is, as discussed in the introductory section,
$y=1+z^{k+1}y^{k+1}$. Upon setting $\lambda=1/y$, we see that this is the same equation. Consequently $D_m$ has an explicit expression
\begin{equation*}
D_m=\alpha y^{-m}+\sum_{i=1}^k \beta_i \rho_i^{m},
\end{equation*}
where the $\rho_i$ are the other roots of the characteristic equation. If one now takes a limit $D_{m+j}/D_m$ for $m\to\infty$ and fixed $j$ the contributions coming from the other roots will disappear, with the result $y^{-j}$. We will use this now:

According to Cramer's rule to solve the linear system of equations, we find
\begin{equation*}
\varphi_i=\frac{D_tz^iD_{h-i}}{D_{h+t+1}},\quad 0\le i\le h
\end{equation*}
and
\begin{equation*}
\varphi_{-i}=\frac{D_hz^iD_{t-i}}{D_{h+t+1}},\quad 0\le i\le t.
\end{equation*}
Since we do not need the upper boundary at level $h$, we push it to infinity:
\begin{equation*}
\lim_{h\to\infty}\frac{D_tz^iD_{h-i}}{D_{h+t+1}}=D_tz^iy^{i+t+1}.
\end{equation*}
We are only interested in the instance $i=0$, with the result $D_ty^{t+1}$. The quantity $\varphi(0)=y_{k;t}$, with lower boundary $-t$ and upper boundary $\infty$ (=no upper restriction) is the generating function of $k_t$-Dyck paths.
So, we obtained for the enumeration of $k_t$-Dyck paths:
\begin{equation*}
y_{k;t}=D_ty^{t+1}=\sum_{0\le \ell \le \frac tk}\binom{t-k\ell}{\ell}z^{(k+1)\ell}(-1)^\ell\cdot y_k^{t+1} .
\end{equation*}
(The summation is over all integers $\ell$ satisfying the inequalities $0\le \ell \le  t/k$.)
In the last expression, we explicitly wrote $y=y_k$ to emphasize the dependency of the generating function on the parameter $k$.
This explains once again that for $0\le t\le k$ we get the simple result $y^{t+1}$. In general, the generating function is given by
\begin{multline*}
y_{k;t}(z)=\sum_{0\le \ell \le \frac tk}\binom{t-k\ell}{\ell}z^{(k+1)\ell}(-1)^\ell\\ \times
\sum_{n\ge0}\frac{t+1}{(k+1)n+t+1}\binom{(k+1)n+t+1}{n}z^{(k+1)n}.
\end{multline*}
The number of $k_t$-Dyck paths of length $(k+1)n$ is then given by
\begin{multline*}
[z^{(k+1)n}]D_ty^{t+1}=\sum_{0\le \ell \le \frac tk}\binom{t-k\ell}{\ell}(-1)^\ell\\ \times
\frac{t+1}{(k+1)(n-\ell)+t+1}\binom{(k+1)(n-\ell)+t+1}{n-\ell}.
\end{multline*}
The polynomial $D_t$ has alternating coefficients, and it is quite likely that there is some sort of an inclusion-exclusion principle underlying.

The polynomial $D_t$ does not have a combinatorial meaning itself, but we may write
\begin{equation*}
\frac{z^t}{D_t}y_{k;t}(z)=z^ty^{t+1}.
\end{equation*}
The generating function $\frac{z^t}{D_t}$ \emph{has} a combinatorial meaning: In the linear system, first
replacing $t$ by 0 and then $h$ and $i$ by $t-1$, leads to $\frac{z^{t-1}}{D_t}$, and it counts 
the $k$-Dyck paths living in the strip $0..h-1$, ending at the highest level $h-1$. The function
$\frac{z^t}{D_t}$ differs only by an extra factor $z$, representing an up-step, touching the level $h$ for the first time. This is exactly the decomposition as given in Fig.~\ref{path5}.

Paths ending on their highest level appear in \cite{yann, beni-master}.

It is interesting to compare $[z^{(k+1)n}]y_{k;t}(z)$ and $[z^{(k+1)n}]y_{k}(z)^t$, by taking their quotient. This is
\begin{equation*}
\sum\limits_{0\le \ell \le \frac tk}\frac{\binom{t-k\ell}{\ell}(-1)^\ell\frac{t+1}{(k+1)(n-\ell)+t+1}\binom{(k+1)(n-\ell)+t+1}{n-\ell}}
{\frac{t+1}{(k+1)n+t+1}\binom{(k+1)n+t+1}{n}}.
\end{equation*}
Taking the limit of this for $n\to\infty$ (only Stirling's formula for factorials is required), we get
\begin{equation*}
\sum_{0\le \ell \le \frac tk}\binom{t-k\ell}{\ell}(-\rho)^\ell, 
\end{equation*}
with 
\begin{equation*}
\rho=\frac{k^k}{(k+1)^{k+1}}.
\end{equation*}
This sum is indeed $=1$ for $0\le t\le k$, but takes smaller values, when $t$ gets larger in relation to $k$. This matches intuition that,
when $t$ is large, the first part of the path that does not  (yet) hit the level $t$ tends to be longer, and the contribution of the rest, which is measured against $[z^{(k+1)n}]y_{k}(z)^t$, tends to be smaller.

\begin{figure}[h]
	\centering
	\includegraphics[width=2.0in]{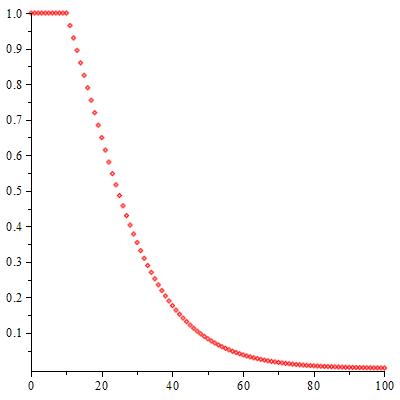}
	\caption{$k=10$, and $t$ is growing}
\end{figure}

%\newpage
\section{Asymptotics}

We want to study the parameter $\mathcal{J}=j$ as in the drawing of Fig.~\ref{path5}. It is given by $j=t+(k+1)\ell$, for some $\ell$. Recall that each path is uniquely decomposed by the $F$-part, until the level $t$ has been reached for the time, and the remainder, the $G$-part.

In order to be able to do explicit calculations, we restrict ourselves to the classical case $k=1$
of Dyck paths. (The case of general $k$, which seems to be less combinatorial, and more analytical, is left as a challenge.)
Then the recursion of second order $D_m=D_{m-1}-z^2D_{m-2}$ admits the solution
\begin{equation*}
D_t(z^2)=D_t(x)=\frac{1-u^{t+1}}{1-u}\frac1{(1+u)^t},
\end{equation*}
with the (classical) substitution $x=\frac{u}{(1+u)^2}$, borrowed from \cite{deBrKnRi72}.

The parameter $\mathcal{J}$ has an automatic contribution of $t$, which we can add later. 
In our setting, we are interested in $\mathcal{J}=t+2\ell$, and we are concentrating on $\ell$. This means that we consider the random variable $\frac{\mathcal{J}-t}2$ and call it the parameter of interest.
The advantage of this procedure is that  now we only have generating functions in $z^2$, for which we write $x$.

The probability generating function of interest is
\begin{equation*}
P(x,w):=\frac{[x^n]\frac1{D_t(xw)}G(x)}{[x^n]C(x)^{t+1}},
\end{equation*}
with
\begin{equation*}
C(x)=\frac{1-\sqrt{1-4x}}{2x}=1+u
\end{equation*}
the generating function of the Catalan numbers, enumerating Dyck paths by half-length: the quantity in the denominator
$[x^n]C(x)^{t+1}$ is the total number of objects (not counting the additional up steps in the beginning); the numerator is
the unique decomposition into the $F$-part and the $G$-part. An additional variable $w$ is used to count the length of the 
$F$-part.

The number of paths of length $2n+t$ with parameter $\mathcal{J}=2s+t$ is then given as $t+2[x^nw^s]P(x,w)$.

We still need to compute $G(x)$, the generating function of paths not going below the line $-t$ and ending on the $x$-axis again.

This can be computed by the linear system as a quotient of the usual determinants: 
Consider paths bounded below by $-t$,  and above by $i$.
\begin{align*}
\lim_{i\to\infty}\frac{D_tD_i}{D_{i+h+1}}&=\lim_{i\to\infty}\frac{1-u^{t+1}}{1-u}\frac{1}{(1+u)^t}
\frac{1-u^{i+1}}{1-u}\frac{1}{(1+u)^i}\\
&\qquad\qquad\qquad\qquad\qquad\qquad\qquad\bigg/
\frac{1-u^{i+t+2}}{1-u}\frac{1}{(1+u)^{i+t+1}}\\
&=\frac{1-u^{t+1}}{1-u}\frac{1}{(1+u)^t}
\frac{(1+u)^{i+t+1}}{(1+u)^i}=\frac{1+u}{1-u}(1-u^{t+1}).
\end{align*}
As a check, we get the product of the two components $F$ and $G$
\begin{equation*}
\frac{1-u}{1-u^{t+1}}(1+u)^t\cdot\frac{1+u}{1-u}(1-u^{t+1})=(1+u)^{t+1}=C(x)^{t+1},
\end{equation*}
as it should.

We study now 
\begin{equation*}
\frac1{D_t(xw)}G(x)
\end{equation*}
using a second variable $w$ to count the parameter of interest.
We compute:
\begin{align*}
\frac{d}{dw}\frac1{D_t(xw)}\Big|_{w=1}&=x\frac{d}{dx}\frac1{D_t(x)}=x\frac{du}{dx}\frac{d}{du}\frac1{D_t(x)}\\
&=-\frac{u}{(1+u)^2}\frac{(1+u)^3}{1-u}\bigg(\frac{(1-u)(1+u)^t}{1-u^{t+1}}\bigg)^2\frac{d}{du}D_t(x)\\
&=-\frac{u(1+u)^{t}}{(1-u^{t+1})^2}\frac{(1+u)(1-u^t)-t(1-u)(1+u^t)}{1-u}\\
\end{align*}
and further
\begin{align*}
\frac{d}{dw}&\frac1{D_t(xw)}\Big|_{w=1}\cdot\frac{1+u}{1-u}(1-u^{t+1})
\\&=-\frac{u(1+u)^{t}}{(1-u^{t+1})^2}\frac{(1+u)(1-u^t)-t(1-u)(1+u^t)}{1-u}\frac{1+u}{1-u}(1-u^{t+1})\\
&=-\frac{u(1+u)^{t+1}}{1-u^{t+1}}\frac{(1+u)(1-u^t)-t(1-u)(1+u^t)}{(1-u)^2}.
\end{align*}
We do not try to simplify this any further, but expand it around $u=1$, which corresponds to the singularity
$x=\frac14$, which is relevant for Dyck-paths. We are in the regime called ``sub-critical'', compare \cite{FS},
where the singular expansion of numerator and denominator is of the same type:
\begin{equation*}
\frac{a_0+a_1\sqrt{1-4x}+\dots}{b_0+b_1\sqrt{1-4x}+\dots},
\end{equation*}
whence the asymptotic expansion of the parameter $\mathcal{J}$ of interest is given by $\frac{a_1}{b_1}$.
One does not need to switch back to the $x$-world, since this quotient can also be obtained via $\frac{a'_1}{b'_1}$, in
\begin{equation*}
\frac{a'_0+a'_1(1-u)+\dots}{b'_0+b'_1(1-u)+\dots}.
\end{equation*}
The reader might recall that $\sqrt{1-4x}=\frac{1-u}{1+u}$, and this is how the singularity $x=\frac14$ translates into $u=1$ and vice versa. A recent application of this (including a proof of a conjectured limiting distribution) is in \cite{Prodinger-hills}.
In our instance, we are led to
\begin{equation*}
\frac{\frac132^tt(t-1)-\frac162^tt(t-1)(t+1)(1-u)+\dots}
{2^{t+1}-2^t(t+1)(1-u)+\dots},
\end{equation*}
and the quotient of interest is
\begin{equation*}
\frac{\frac162^tt(t-1)(t+1)}{2^t(t+1)}=\frac{t(t-1)}{6}.
\end{equation*}
Multiplying this by 2, in order to switch from half-length to length, and adding the fixed contribution $t$ leads to the average value of $\mathcal{J}$:
\begin{align*}
t+2\frac{t(t-1)}{6}=\frac{t(t+2)}{3}.
\end{align*}

In this subcritical regime it is also relatively easy to determine the discrete limiting distribution. 
We start from
\begin{equation*}
\frac{\frac1{D_t(xw)}G(x)}{C(x)^{t+1}};
\end{equation*}
the quotient $G(x)/C(x)^{t+1}$ does not depend on the parameter, and can thus be replaced by its limit at the singularity $x=\frac14$, or, easier at $u=1$: (This simplification stems from the fact that the limit of a product is the product of limits, and, more generally, this holds for the local expansions as well.)
\begin{equation*}
\lim_{u\to1}\frac{G(x)}{C(x)^{t+1}}=\frac{2(t+1)}{2^{t+1}}=\frac{t+1}{2^t}.
\end{equation*}
Hence, the discrete limiting distribution is given by the (probability) generating function
\begin{equation*}
p_t(u)=\frac{t+1}{2^t}\frac{1-u}{1-u^{t+1}}(1+u)^t.
\end{equation*}
One can even read off coefficients explicitly: (The contour is again a small circle in the $x$-plane, and since $x\approx u$ for small $x$, the curve also winds around the origin exactly once in the $u$-plane.)
\begin{align*}
[x^m]&p_t(u)=\frac{t+1}{2^t}\frac1{2\pi i}\oint\frac{dx}{x^{m+1}}p_t(u)\\
&=\frac{t+1}{2^t}\frac1{2\pi i}\oint\frac{du(1-u)(1+u)^{2m+2}}{u^{m+1}}\frac{1-u}{1-u^{t+1}}(1+u)^t\\
&=\frac{t+1}{2^t}[u^m]\frac{(1-u)^2}{1-u^{t+1} }(1+u)^{2m-1+t}\\
&=\frac{t+1}{2^t}\sum_{\lambda\ge0}\bigg[\binom{2m-1+t}{m-\lambda(t+1)}\\
&\qquad\qquad\qquad-2\binom{2m-1+t}{m-1-\lambda(t+1)}+\binom{2m-1+t}{m-2-\lambda(t+1)}\bigg].
\end{align*}
This quantity can be interpreted as the probability that in an (almost) infinitely long path the parameter
$\frac{\mathcal{J}-t}{2}$ has value $m$.

\begin{figure}[h]
	\centering
	\includegraphics[width=1.45in]{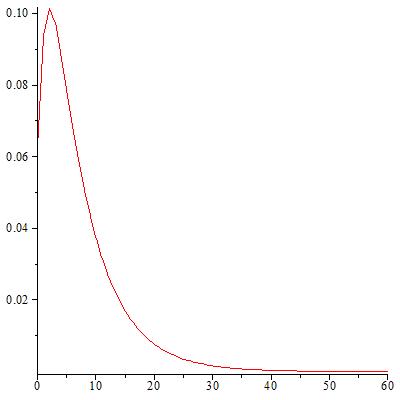}
	\includegraphics[width=1.45in]{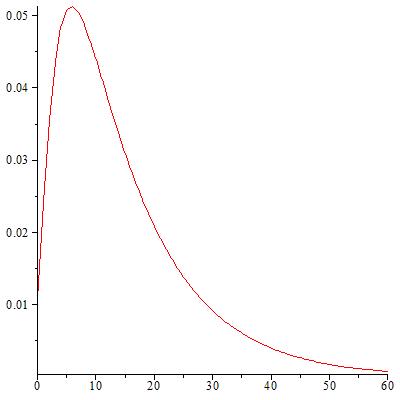}
	\includegraphics[width=1.45in]{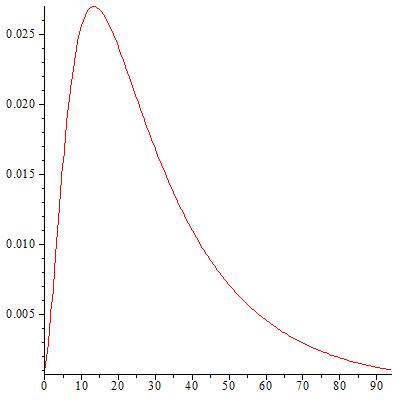}
	\caption{ Limiting distribution: $t=7$, $t=10$,  $t=14$ }
\end{figure}

\textbf{Acknowledgments.} I thank Nancy Gu for valuable discussions. Most of this work was developed while I visited the Technical University in Graz, Austria. I thank my host Peter Grabner for his hospitality. Furthermore, I thank two referees who pointed out several formulations and statements that lacked in clarity. In this revised version I provided much more detailed explanations and fixed some inaccuracies.

\newpage

\bibliographystyle{plain}
%\bibliography{k-t-Dyck}

\end{document}